\documentclass[12pt,amssymb,epsf]{article}
\usepackage{epsfig}
\usepackage{emlines}
\usepackage{epic,eepic}
\usepackage{amsfonts,epsfig,latexsym}
\usepackage{emlines}

\title{ \bf Clusters of Cycles}

\author{Michel DEZA \thanks{first author acknowledges financial support of Mathematical Institute at Ohio State University}\\
     CNRS and Ecole Normale Sup\'erieure, Paris, France \\
\and   Mikhail SHTOGRIN \thanks{second author acknowledges financial support 
of the Russian Foundation for Fundamental Research (grant 99-01-00010)}\\
Steklov Mathematical Institute, Moscow, Russia} 
\date{\today}

\begin{document}

\maketitle

\begin{abstract}

A { \it cluster of cycles} (or {\it $(r,q)$-polycycle}) is a simple planar 2--connected finite or countable graph
$G$ of girth $r$ and maximal vertex-degree $q$, which admits an
{\it $(r,q)$-polycyclic realization} $P(G)$ on the plane. An 
$(r,q)$-polycyclic realization is determined by the following properties:
(i) all interior vertices are of degree $q$,
(ii) all interior faces (denote their number by $p_r$) are combinatorial $r$-gons,
(iii) all vertices, edges and interior faces form a cell-complex.

An example of $(r,q)$-polycycle is the skeleton of $(r^q)$, i.e. of the
$q$-valent partition of the sphere, Euclidean plane or
hyperbolic plane by regular $r$-gons. Call
{\it spheric} pairs $(r,q)=(3,3),(4,3),(3,4),(5,3),(3,5)$. Only for those
five pairs, 
$P((r^q))$ is $(r^q)$ without exterior face; otherwise, $P((r^q))=(r^q)$.

Here we give a compact survey of results on $(r,q)$-polycycles.

We start with the following general results for any $(r,q)$-polycycle $G$:
(i) $P(G)$ is unique, except of (easy)
case when $G$ is the skeleton of one of 5 Platonic polyhedra;
(ii) $P(G)$ admits a cell-homomorphism $f$ into $(r^q)$;
(iii) a polynomial criterion to decide if given finite graph is a polycycle, is presented.

Call a polycycle {\it proper} if it is a partial subgraph of
$(r^q)$ and a {\it helicene}, otherwise.
In [\ref{Ha}] all proper spheric
polycycles are given.
An $(r,q)$-helicene exists if and only if 
$ p_r > (q-2)(r-1)$ and $(r,q) \neq (3,3)$. We list the $(4,3)$-, $(3,4)$-helicenes and
the number of $(5,3)$-, $(3,5)$-helicenes for first interesting $p_r$.
Any outerplanar $(r,q)$-polycycle $G$ is a proper $(r,2q-2)$-polycycle
and its projection $f(P(G))$ into $(r^{2q-2})$ is convex.
Any outerplanar $(3,q)$-polycycle $G$ is a proper $(3,q+2)$-polycycle.

The symmetry group $Aut(G)$ (equal to $Aut(P(G))$, except of Platonic case) of an $(r,q)$-polycycle $G$ is a subgroup of
$Aut((r^q))$ if it is proper and an extension of $Aut(f(P(G)))$, otherwise. 
$Aut(G)$ consists only of rotations and mirrors if $G$ is finite; so its order
divides one of the numbers $2r$, $4$ or $2q$. Almost all polycycles $G$ have trivial $AutG$.

Call a polycycle $G$ {\it isotoxal} (or {\it isogonal}, or
{\it isohedral}) if $AutG$ is transitive on edges (or vertices, or interior
faces); use notation IT (or IG, or IH), for short. Only $r$-gons and non-spheric $(r^q)$ are isotoxal. Let $T^*(l,m,n)$ denote Coxeter's triangle group of a 
triangle on $S^2$, $E^2$ or $H^2$ with angles $ \frac{ \pi}{l}$, $\frac{ \pi}{m}$, $ \frac{ \pi}{n}$ and
let $T(l,m,n)$ denote its subgroup of index 2, excluding motions of 2nd kind.
We list all IG- or IH-polycycles for spheric $(r,q)$ and construct many examples
of IH-polycycles for general case (with $AutG$ being above two groups for some parameters,
including strip and modular groups). Any IG-, but not IT-polycycle is infinite,
outerplanar and with same vertex-degree; we present two IG-, but not
IH-polycycles with $(r,q)=(3,5),(4,4)$ and $AutG=T(2,3, \infty) \sim PSL(2,Z),
T^*(2,4, \infty)$. Any IH-polycycle has the same number of boundary
edges for each its $r$-gon. For any $r \ge 5$ there exists a continuum of
{\it quasi}-IH-polycycles, i.e. not isohedral, but all $r$-gons have the same 1-corona.

On two notions of extremal polycycles:

(i) We found, for the spheric $(r,q)$ the maximal number $n_{int}$ of interior points for a
$(r,q)$-polycycle with given $p_r$; in general case,
$ \frac{p_r}{q} \le n_{int} < \frac{rp_r}{q}$ if any $r$-gon contains an interior point.

(ii) All {\it non-extendible} $(r,q)$-polycycles (i.e. not a proper subgraphs
of another $(r,q)$-polycycle) are $(r^q)$, 4 special ones,
(possibly, but we conjecture their nonexistence) some other finite $(3,5)$-polycycles, and, for any
$(r,q) \neq (3,3),(3,4),(4,3)$, a continuum of infinite ones.

On isometric embedding of polycycles into hypercubes
$Q_m$, half-hypercubes $ \frac{1}{2}Q_m$ and, if infinite, into cubic lattices
$Z_m$, $ \frac{1}{2}Z_m$: for $(r,q) \neq (5,3),(3,5)$,
there are exactly 3 non-embeddable polycycles (including $(4^3)-e$, $(3^4)-e$); all non-embeddable $(5,3)$-polycycles
are characterized by two forbidden sub-polycycles with $p_5=6$.
\end{abstract}



\pagestyle{plain}

\noindent

\section{Definition and examples  }

A { \it cluster of cycles} (a { \it polycycle}, for short, or {\it $(r,q)$-polycycle}) is a simple planar 2-(vertex)-connected finite or countable graph
$G$ of girth $r$ and maximal vertex-degree $q$, which admits an
{\it $(r,q)$-polycyclic realization} $P(G)$ on the plane.

An $(r,q)$-polycyclic realization is determined by the following properties:

(i) all interior vertices are of degree $q$,

(ii) all interior faces are  (combinatorial) $r$-gons,

(iii) all vertices, edges and interior faces form a cell-complex (i.e. the
intersection of any two faces is edge, vertex or $\emptyset$).

\vspace{2mm}
One can show that (iii) follows from (i), (ii), while neither (i), (iii) imply (ii), nor (ii), (iii)
imply (i).

\vspace{2mm}
For example, $(3,q)$- and $(4,q)$-polycycles
are just simplicial and cubical complexes of dimension 2. 

The main example of $(r,q)$-polycycle is the skeleton of $(r^q)$, i.e. of the
$q$-valent partition of the sphere $S^2$, Euclidean plane $R^2$ or
hyperbolic plane $H^2$ by regular $r$-gons. For  $(r,q) = (3,3), (4,3), (3,4),
(5,3) ,(3,5)$, the unique $(r,q)$-polycycle is, respectively, Platonic tetrahedron, cube,
octahedron, dodecahedron, icosahedron on $S^2$, {\it but with excluded
exterior face}; for $(r,q) = (6,3), (3,6),
(4,4)$ it is regular partition $(6^3), (3^6), (4^4)$ of $R^2$; all others
$(r^q)$ are regular partitions of $H^2$.

Call a polycycle {\it proper} if it is a partial subgraph of (the skeleton of)
the regular partition $(r^q)$. Call a proper $(r,q)$-polycycle
{\it induced} (moreover, {\it isometric}) if it is induced (moreover,
isometric) subgraph of $(r^q)$.

 Call an $(r,q)$-polycycle {\it spheric}, {\it Euclidean} or 
{\it hyperbolic}, if $(r^q)$ is the
regular partition of $S^2$, $R^2$ or $H^2$, respectively. (One can also use 
terms {\it elliptic}, {\it parabolic} or {\it hyperbolic}, since 
$rq<2(r+q)$, $rq=2(r+q)$ or $rq>2(r+q)$, respectively.) There is 
a literature (see, for example, Section 9.4 of [\ref{GS}] and [\ref{BGOR}])
about proper Euclidean polycycles ({\it polyhexes}, {\it polyamonds},
{\it polyominoes} for $(6^3), (3^6), (4^4)$, respectively); the terms come
from usual terms {\it hexagon}, {\it diamond}, {\it domino}, where the last 
two correspond to the case $p_3$, $p_4=2$.  Polyominoes
were considered first by Conway, Penrose, Colomb as tilers (of $R^2$ etc.; 
see, for example, [\ref{CL}]) and
in the games; later they were used for enumeration in Physics and Statistical
Mechanics. Polyhexes are used widely (see, for example, [\ref{Di}], [\ref{Ba}]) in
Organic Chemistry: they represent {\it completely condensed PAH (polycyclic
aromatic hydrocarbons)} $C_nH_m$ with $n$ vertices (atoms of the carbon $C$),
including $m$ vertices of degree 2, where atoms of the hydrogen $H$ are
adjoined. All 39 proper $(5,3)$-polycycles were found in [\ref{CC}] in chemical
context, but already in [\ref{Ha}] were given all 3, 6, 9, 39, 263 proper spheric
$(r,q)$-polycycles for  $(r,q)=(3,3), (4,3), (3,4), (5,3),(3,5)$,
respectively. 

A general theory of polycycles is considered in [\ref{DS3}]-[\ref{DS8}], 
[\ref{S1}]-[\ref{S2}]. Namely, proofs can be found as follows: Theorem 1 in 
([\ref{S1}]-[\ref{S2}]), Theorems 2-3 in ([\ref{DS3}], [\ref{DS5}], 
[\ref{DS6}]), Theorems 4-6 in [\ref{DS7}], Theorems 7-9 in 
[\ref{DS8}], Theorem 10 in ([\ref{DS4}], [\ref{DS7}]).
 
\section{Criterion and unicity}

\vspace{2mm}
\noindent
{\bf Theorem 1.} {\it Let $G$ be any {\it finite} connected
graph of girth $r$, different from the skeleton of $(3^3), (3^4), (4^3),
(3^5), (5^3)$; let $v, e, f$ be its number of vertices, edges and $r$-cycles,
respectively. Then $G$ is $(r,q)$-polycycle if and only if holds:

(i) any edge belongs to one or two $r$-cycles of the graph $G$,

(ii) all edges, belonging to exactly one $r$-cycle of $G$, form a
 simple cycle,

(iii) the intersection of any two different $r$-cycles is an edge, a vertex or
$\emptyset$,

(iv) $v-e+f=1$,

(v) all $r$-cycles with common vertex can be organized into a sequence, such 
that any two neighbors have the common edge, containing the common
vertex, and this sequence has at most $q$ members with equality if and only if
it is closed (i.e. the sequence form a cycle).}

\vspace{2mm}
It is clear that: (i) implies that $G$ is 2-connected, (v) implies that any
interior vertex of $G$ has degree $q$ and that for $q=3$ the condition (v) is
implied by (i)-(iv) with exclusion in (iii) of the case of intersection in a
vertex.

\vspace{2mm}
\noindent
{\bf Theorem 2.} {\it Let $G$ be an $(r,q)$-polycycle. Then
 
(i) if $G$ is the skeleton of one of 5 Platonic polyhedra, then the number
of $(r,q)$-polycyclic realizations of $G$ is equal to the number of faces of
the Platonic polyhedron and all those realizations are isomorphic,

(ii) any other polycycle $G$ has unique $(r,q)$-polycyclic realization and the number of its interior faces (which are all should be $r$-gons) is the number of $r$-cycles of $G$.}

\vspace{2mm}
The $(r,q)$-polycyclic realization $P(G)$ of a non-Platonic $(r,q)$-polycycle
$G$ is, in general, not unique {\it planar} realization of polycycle $G$.

\section{Proper polycycles versus helicenes}

\vspace{2mm}
First, we list all $(3,3)-,(4,3)-,(3,4)$-polycycles. 
Clearly, all $(3,3)$-polycycles are $(3^3)$, $(3^3)$-$v$ and $(3^3)$-$e$ (i.e.
a vertex with incident edges, or an edge is deleted); the last one is not induced.

We denote by $P_n$ a path with $n$ vertices and by $P_N$, $P_Z$ infinite paths
in one or both directions. All $(4,3)$-polycycles are: 
$(4^3)$, $(4^3)$-$v$, $(4^3)$-$e$, $P_2 \times P_n$ for any natural $n$ 
and two infinite ones: $P_2 \times P_N$, $P_2 \times P_Z$.
Only  $(4^3)$, $(4^3)$-$v$, $P_2 \times P_2$,
$P_2 \times P_3$, $P_2 \times P_4$, $(4^3)$-$e$ are proper; the last two are
not induced. 

The number of $(3,4)$-polycycles is also countable, including two
infinite ones (9 of $(3,4)$-polycycles are proper and 5 of proper ones are induced). Namely, all $(3,4)$-polycycles are: proper ones $(3^4)$, $(3^4)$-$v$,
$(3^4)$-$e$, $(3^4)$-$P_3$, $(3^4)$-$C_3$, $G_n$ ($1 \le n \le4$) and unproper
ones $G_n$ ($n \ge 5$ and $n = N, Z$) and the {\it vertex-split $(3^4)$}, defined below. Here
$G_n = A_{ \frac{n}{2}}$ if $n$ even, $= B_{ \frac{n+1}{2}}$ if $n$ odd,
$= A_n$ if $n = N, Z$, where $A_m$  is  $(4,3)$-polycycle  $P_2 \times P_m$
with parallel diagonals, added one on each square, and $B_m$ is
$A_m$ without a vertex of degree 2. 

For all other $(r,q)$, 
there is a continuum of $(r,q)$-polycycles and the number of finite ones 
among them is countable.

Call {\it vertex-split $(3^4)$}, a $(3,4)$-polycycle, 
coming from $(3^4)$ as follows: let $K_{x, \{ a,b,c,d \}}$ be induced 4-wheel 
in $(3^4)$, then replace the edges $(x,a), (x,b)$ of $(3^4)$ by the edges 
$(x',a), (x',b)$, where $x'$ is a new vertex of degree 2. (Curiously, 
this plane graph is the logo of the HSBC, Hong Kong and Shanghai Banking 
Corporation.) 

Call {\it vertex-split $(3^5)$}, a $(3,5)$-polycycle, 
coming from $(3^5)$ as follows: let $K_{x, \{ a,b,c,d,e \}}$ be induced
5-wheel in $(3^5)$, then replace the edges $(x,a), (x,b)$ of $(3^5)$ by the 
edges $(x',a), (x',b)$ where $x'$ is a new vertex.
 
\subsection{Cell-homomorphism into $(r^q)$ and helicenes}

\vspace{2mm}
\noindent
{\bf Theorem 3.} {\it Any
$(r,q)$-polycyclic realization $P(G)$ admits a cell-homomorphism into $(r^q)$ and
such homomorphism is defined uniquely by a {\it flag} (i.e.
incident vertex, edge and interior face of $P(G)$) and its image.}

\vspace{2mm}
Clearly, the above homomorphism is an isomorphism if and only if $G$ is a proper
polycycle (i.e. if the map $P(G) \to (r^q)$ is topologic: there is no pair of 
vertices or of edges, having the same image). In view of Theorem 3 any 
unproper $(r,q)$-polycycle is called {\it $(r,q)$-helicene}.

\vspace{2mm}
It is easy to check that $(r,q)$-helicenes exist if and only if $(r,q) \neq (3,3)$ 
and $p_r \ge (q-2)(r-1)+1$ with equality only for the helicene being a belt
of $r$-gons, going arond an $r$-gon. All
$(4,3)$-helicenes are two infinite ones: $P_2 \times P_N$, $P_2 \times P_Z$
and $P_2 \times P_n$ for any $n \ge 5$; given above full description of all
$(3,4)$-polycycles, permit also to list all $(3,4)$-helicenes. We counted
that the number of $(5,3)$-helicenes is $1,7,29$ for $p_5=5,6,7$ and  
the number of $(3,5)$-helicenes is $1,4,20,74$ for $p_3=7,8,9,10$.

A natural
parameter to measure an $(r,q)$-helicene, will be the degree of the 
corresponding homomorphism into $(r^q)$ (on vertices, edges and faces). 
For $q \ge 4$, helicenes appear with vertices, but not edges,
having same homomorphic image. The vertex-split $(3^4)$ is unique
such maximal helicene for $(r, q) = (3, 4)$ ($x, x'$ are such vertices). There is
a finite number of such helicenes for $(r, q) = (3, 5)$; one of them is the
vertex-split $(3^5)$.

\subsection{Cell-complexes $P(G)$, $K(G)$, $X(G)$ and the curvature}

Denote by $K(G)$ the abstract 2-dimensional polyhedron with a metric, such
that all $r$-gons of $P(G)$ became planar Euclidean regular $r$-gons; clearly,
there is combinatorial cell-isomorphism between $K(G)$ and $P(G)$. The map
$K(G) \to K(r^q)$ is a geometric realization of the combinatorial map
$P(G) \to (r^q)$, such that the homomorphism is locally-homeomorphic, i.e.
it is continuous cell-map with isomorphic $\epsilon$-neighborhood for
sufficiently small $ \epsilon$ ({\it homeomorphic} means {\it isometric} and
{\it global} homeomorphism means {\it isomorphism}).

The gaussian curvature of a point in $K(G)$ is 
$2 \pi - q \frac{r-2}{r} \pi $ in each interior vertex (since in each 
interior vertex meet $q$ angles of regular Euclidean $r$-gons) and $0$ in
any other point (since $K(G)$ is a disc glued from Euclidean $r$-gons). So
the global curvature of $K(G)$ is the sum of its curvatures in interior vertices:

$n_{int}(2 \pi - q \frac{r-2}{r} \pi)=n_{int} \frac{\pi}{r} (2(r+q)-rq)$.

For example, the curvature of $G$ is $3 \pi$ for $G= (3^5), (5^3)$,
is $2 \pi$ for $G=(4^3), (3^4)$ and  $\pi$ for $G=(3^3)$ 
(while the curvature of the sphere $S^2$ is $4 \pi$).

\vspace{2mm}
Denote by $X(G)$ the metric space of {\it constant} curvature, obtained from $P(G)$
by introducing a metric on it, which is locally spheric, locally Euclidean or locally
hyperbolic, if $(r^q)$ is a regular partition of $S^2$, $R^2$ or $H^2$, respectively.
$X(G)$ has also cell structure, glued from, in general, non-Euclidean faces, but here
we consider it as an abstract cell-complex. Clearly, $X(G)=K(G)$ for $(r,q)=(4,4), (6,3), (3,6)$.
In general, both complexes are homeomorphic as manifolds and have the same 
curvature, but it is
non-zero only on interior vertices of $K(G)$ and it is constant on all 
points of $X(G)$. There is cell-isomorphism amongst complexes
$X(G)$, $K(G)$, $P(G)$ and each of them admits cell-homomorphism on corresponding complex
of $(r^q)$. $X(G)$ is also simply-connected 2-dimensional manifold, which is homeomorphic to
a disc if $G$ is finite and non-compact, otherwise. The manifold $X(G)$ has no boundary only if
$G$ is the skeleton of partition $(r^q)$ of Euclidean or hyperbolic plane.
All faces of the complex $X(G)$ are regular $r$-gons
with angles $\frac{2 \pi}{q}$, while the faces of complex $K(G)$ are regular 
Euclidean $r$-gons with angles $\frac{(r-2) \pi}{r}$.

\subsection{Outerplanar polycycles}

Call a polycycle {\it outerplanar} if it has no interior points, i.e.
$n_{int}=0$; clearly, it is a $(r,q')$-polycycle for any $q'$ not less than
the maximal degree of vertices. The following Theorem show that, in a way,
outerplanar polycycles are close to proper polycycles.

\vspace{2mm}
\noindent
{\bf Theorem 4.} 

{\it (i) any outerplanar $(r,q)$-polycycle $G$ is a proper $(r,2q-2)$-polycycle 
and its projection $f(P(G))$ into $(r^{2q-2})$ is convex (on $S^2$, $R^2$ or $H^2$),
 
(ii) any outerplanar $(3,q)$-polycycle is a proper $(3,q+2)$-polycycle.}

\vspace{2mm}
The proof uses the fact that the projection on $(r^q)$ of polycyclic realization of the graph, being simply-connected, is convex if and only if all boundary angles
are less than $\pi$ (the boundary will be a union of convex polygons).

Remark that already for $p_3=7$  there are outerplanar $(3,4)$- and
$(3,5)$-polycycles, which remain unproper in  $(3^5)$, $(3^6)$, respectively.
A {\it fan} of $(q-1)$ $r$-gons with $q$-valent common (boundary) vertex, is an example of
outerplanar $(r,q)$-polycycle, which is a proper {\it non-convex}
$(r,2q-3)$-polycycle.

\subsection{Proper and reciprocal polycycles}

For a proper polycycle we are interested when it is induced 
(or, moreover, isometric) subgraph of $(r^q)$. For $(r,q)=(3,3),(4,3),(3,4)$ 
any induced $(r,q)$-polycycle for $(r,q)=(3,3),(4,3),(3,4)$ is isometric, but, for
example, the path of 3 pentagons is induced non-isometric $(5,3)$-polycycle.
Any isometric polycycle is embeddable (see Section 7 below), but already
for $p_5=6$ there exists a non-embeddable induced $(5,3)$-polycycle.

Other possible property of a proper $(r,q)$-polycycle is being convex in $(r^q)$
(see Theorem 4 and remark after Lemma below). 
Consider now the notion of reciprocity, defined for some proper polycycles.

Let $P$ be a proper bounded $(r,q)$-polycycle. Consider the union of all 
($r$-gonal) faces of $(r^q)$ outside of $P$.
Easy to see that this union will be an $(r,q)$-polycycle (and call it then
{\it reciprocal to} $P$) if and only if, either $P$ is spheric, or $P$ is
infinite and has connected boundary. Call a polycycle {\it self-reciprocal}
if it admits the reciprocal polycycle and is isomorphic to it.

All self-reciprocal $(r,q)$-polycycles with $(r,q)=(3,3),(4,3),(3,4),(5,3)$
are  $(3^3)$-$e$, $(4^3)$-$v$, $P_2\times P_3$, $(3^4)$-$v$, $(3^4)$-$C_3$,
$(3^4)$-$2K_2$ and 9 (out of 11) $(5,3)$-polycycles with $p_5=6$, including
6 chiral ones. An example of self-reciprocal $(3,q)$-polycycle for any $q \ge 3$, is a $(3,q)$-polycycle on one of two shores of {\it zigzag path}, cutting
$(3^q)$ in two isomorphic halves; it includes $(3^3)$-$e$ and $(3^4)$-$C_3$
and infinite for $q \ge 6$.

\section{Symmetries of polycycles}

The symmetry group $AutG$ of an $(r,q)$-polycycle $G$ is a subgroup of
$Aut(r^q)$ if it is proper and an extension of $Aut(HomG)$,
otherwise; here $HomG$ denotes the cell-homomorphism projection of $P(G)$
into $(r^q)$. We have $AutG=AutP(G)$, except of the case of $G$ being one of five Platonic $(r^q)$. 
If an $(r,q)$-polycycle $G$ is finite and $P(G)$ has a fixed point inside it, then
$Aut(G)$ consists only of rotations and mirrors around of this point. So its order
divides $2r$, $4$ or $2q$, depending on what $AutG$
fixes: the center of an $r$-gon, the center of an edge or a vertex (the corresponding
groups are $D_{rh}$, $D_{2h}$, $D_{qh}$). (Above $AutG$ is given, for finite 
polycycle $G$, as a space group, i.e. we discard plane mirror.) None of $(3,3)$-, $(3,4)$-, 
$(4,3)$-polycycles, but almost all $(r,q)$-polycycles for any other $(r,q)$,
have trivial $AutG$. The number of {\it chiral} (i.e. with AutG containing no mirrors)
proper $(5,3)$-, $(3,5)$-polycycles is 12, 208 (amongst, respectively, all 39,
263).

Call a polycycle $G$ {\it isotoxal} (or {\it isogonal}, or
{\it isohedral}) if $AutG$ is transitive on edges (or vertices, or interior
faces); use notation IT- (or IG-, or IH-), for short.

Let $T^*(l,m,n)$ denote Coxeter's {\it triangle group} of a fundamental triangle
with angles $ \frac{ \pi}{l}$, $\frac{ \pi}{m}$,
$ \frac{ \pi}{n}$. Let $T(l,m,n)$ denote its subgroup of index 2, excluding
motions of 2nd kind (i.e. those, changing orientation); see, for example [\ref{Ma}], pages 81,90,107,176,183.
Now, $T^*(2,2, \infty )=pmm2$, $T(2,2, \infty )=p112 \approx pm11 \approx pma2$. 
(Remark that $p1m1$ also has index 2 in $T^*(2,2, \infty )$, but it is not 
isomorphic to $T(2,2, \infty )$.) On the other hand, 
$T(2,3, \infty ) \approx PSL(2,Z)$ (the modular group) and 
$T^*(2,3, \infty )  \approx SL(2,Z)$.
For all but one (the last on Figure 4) known families of infinite IG- or 
IH-polycycles, $AutG$, if it is not a strip group, is one of above two groups
$T(l,m,n)$, $T^*(l,m,n)$ for some parameters. For all known such polycycles with strip group $AutG$
(see Figure 3), this group is isomorphic to one of $T^*(2,2, \infty )$,
$T(2,2, \infty )$.

Only $r$-gons and non-Platonic $(r^q)$
are isotoxal; their respective symmetry groups are $D_{rh}$ and $T^*(l,m,n)$.
$Aut(r^q)=D_{rh}$ in five Platonic cases; none is IT-, IG- or IH-polycycle, 
except of isohedral $(3^3)$. 

We conjecture that all polycycles, which are isogonal and isohedral, but not
isotoxal, are the infinite $(3,4)$-polycycle from 3rd column in Figure 1 and
$(2k,3)$-{\it cactuses} for any $k \ge 2$ (with $AutP=T^*(2,k, \infty ))$ and
we checked this conjecture for spheric $(r,q)$. In fact, 
only other IG-, but not IT-polycycle in Figure 1 of all spheric IH-polycycles,
is $P_2 \times P_Z$, i.e. the $(4,3)$-cactus.
The cactuses are infinite polycycles obtained by the procedure, which is clear
from Figures 2 and 4. The $(2k,3)$-cactuses for any $k \ge 3$ correspond to the
case $a=0$ of the first family in Figure 4; they are only isogonal polycycles
in Figure 4.

\vspace{2mm}
\noindent
{\bf Theorem 5.}
 
{(i)\it Any IG-, but not IT-polycycle is infinite,
outerplanar and with the same vertex-degree;

(ii) there exist two IG-, but not
IH-polycycles with $(r,q)=(3,5),(4,4)$ (see Figure 2; their groups are 
$T(2,3, \infty)$, $T^*(2,4, \infty)$) and this $(3,5)$--polycycle is unique 
such spheric polycycle.}

\noindent
{\bf Theorem 6.}
{\it Let $P$ be an isohedral $(r,q)$-polycycle. Then

(i) $P$ has the same number $t$ of non-boundary
edges for each its $r$-gon;

(ii) if $t=0, r$ or $1$, then $P$ is $r$-gon, $(r^q)$ (with $(r,q) \neq (4,3),(3,4),(5,3),(3,5)$) or a pair of adjacent $r$-gons (and $AutP=D_{2h}$);

(iii) if $t=2$, then $P$ is either a {\it star} of $q$ $r$-gons with one 
common interior vertex, or an infinite outerplanar polycycle;

(iii') there exist exactly two infinite isohedral $(3,q)$-polycycles: both 
infinite polycycles from the 5th column in Figure 1 (both with 
$AutP=pma2 \approx T(2,2, \infty )$);

(iv) if $1 \le t \le r-3$, then $P$ is infinite; for any $r$ there exists
isohedral $(r,4)$-polycycle with $t=r-1$(take an $r$-gon with right angles and $AutP$, generated by mirrors from all but one its sides); 
  
(v) all spheric isohedral polycycles are (see Figure 1): 11 finite ones (see 
(ii) above) and 8 infinite ones ($P_2 \times P_Z$, six its decorations (all with strip groups $AutP$) and 
the $(5,3)$-cactus with $AutP=T^*(2,3, \infty )$;

(vi) 8 families of isohedral decorations of $P_2 \times P_Z$ are given in Figure 3; 9 families of isohedral decorations of $(r,q)$-cactuses are given in Figure 4.}

\vspace{1mm}
Remark that 1st, 2nd decorations in Figure 3 are the case $k=2$ of, 
respectively, 1st, 4th cactuses on Figure 4. Amongst the 8 decorations in Figure 3,
only the case $a=0$ of the 4th one is isogonal.
 
The group of the last family in Figure 4 is in 1-1-correspondence with
$T(2,2, \infty )$, but different from it. It is the product of
$T^*(\infty, \infty, \infty )$ and the group (of order 3) of congruence
of the triangle with all vertices in infinity.

\vspace{2mm}
For any $r \ge 5$ there exists a continuum of {\it quasi}-IH-polycycles, i.e. 
not isohedral, but all $r$-gons have the same 1-corona.
In fact, let  $T$  be an infinite, in both directions, path of regular 
$r$-gons, such that for any of them, the edges of adjacency to its neighbors 
are at distance $ \lfloor \frac {r-3}{2} \rfloor $ and the sequence of (one of
two possible) choices of joining each new $r$-gon, is aperiodic and different
from its reversal. There is a continuum of such $T$ for any $r \ge 5$.
Any $T$ is quasi-isohedral and its group of automorphism is trivial. Also,
$T$ is embeddable (see Section 7 below); it is unproper if $r = 5, 6$ and 
isometric if $r \ge 7$.  


\newpage
\vspace{6mm}

\begin{center}
\input tab4.pic
\end{center}
\begin{center}
{Figure 1: all isohedral spherical polycycles (only                            $r$-gons and two infinite ones with $r, q$ different from $5$, are isogonal)}
\end{center}

\newpage
\vspace{6mm}

\begin{center}
\input tab5.pic
\end{center}
\begin{center}
{Figure 2: examples of isogonal but not isohedral polycycles }
\end{center}

\vspace{6mm}
\begin{center}
\input tab6.pic
\end{center}
\begin{center}
{Figure 3: examples of isohedral polycycles with strip groups}
\end{center}

\newpage
\vspace{6mm}

\begin{center}
\input tab7.pic
\end{center}
\begin{center}
{Figure 4: examples of families of isohedral cactuses}
\end{center}

\newpage
\section{Two extremal properties}

\subsection{Maximal number of interior points}

Recall that $p_r(P)$, $n_{int}(P)$ denote the number of interior faces and interior vertices
of given finite $(r,q)$-polycycle $P$. Call $\frac{n_{int}(P)}{p_r(P)}$ the {\it density} of $P$ and denote by $n(x)$ the maximum of $n_{int}(P)$ over all
$(r,q)$-polycycles $P$ with $p_r(P)=x$; call {\it extremal} any
$(r,q)$-polycycle $P$ with $n_{int}(P)=n(x)$. So, extremal
polycycles represent opposite case to outerplanar ones (see Section 3.3. 
above)), with the same $p_r$. Clearly, an 
extremal polycycle also maximises the number $e_{int}$ of 
non-boundary edges and minimizes the number $Per$ (for 
{\it perimeter})
of boundary edges (or boundary points), as well as the number $n$ of all
vertices and the number $e$ of all edges. In fact, Euler formula

$(n_{int}+Per)-(e_{int}+Per)+(p_r+1)=2$ 

and equality $rp_r=2e_{int}+Per$, imply

$n_{int}=e_{int}-p_r+1=-\frac{Per}{2} +\frac{p_{r}(r-2)}{2}+1=-e+p_r(r-1)+1=-n+p_r(r-2)+2$.

For $(5,3)$-polycycles with $x \le 11$, $n(x)$ was found in [\ref{CC}; all 
extremal $(5,3)$-polycycles turn out to be proper and unique. Moreover, the
$(5,3)$-polycycles, which are reciprocal to any such extremal one, turn out to be also extremal. [\ref{CC}] asked 
about $n(x)$ for $x \ge 12$; this Section answer this question for any $x$ and
for all spheric $(p,q)$.

First, consider three trivial cases of $(r,q)$. All pairs 
$(p_r,n_{int})$ for $(p,q)=(3,3)$ are $(1,0)$, $(2,0)$, $(3,1)$; for 
$(p,q)=(4,3): (n,0)$ for all $n \ge 1$, $(3,1)$, $(4,2)$, 
$(5,4)$ and for $(p,q)=(3,4): (n,0)$ for all $n \ge 1$,
$(4,1)$, $(5,1)$, $(6,1)$, $(6,2)$, $(7,3)$.

Call {\it kernel} of a polycycle, the cell-complex of vertices, edges and 
faces of the polycycle, which are not incident with its boundary. Call a polycycle
{\it elementary} if it is a $r$-gon or if it has non-empty connected kernel, such that
the deletion of any face from the kernel will diminish it (i.e. any face of the
polycycle is incident to its kernel).

\vspace{2mm}
\noindent
{\bf Lemma.}\ {\it 

(i) Any $(r,q)$-polycycle with $(r,q)=(3,3),(3,4),(4,3),(3,5),
(5,3)$ is a union of elementary polycycles without common faces, 

(ii) for $(r,q)=(6,3),(4,4),(3,6)$ and for the case
$r \ge 7, q \ge 3$, there is a continuum of connected (components of) kernels of $(r,q)$-polycycles.}

In above Lemma (ii) and for $r \ge 7, q \ge 4$, one can find,
moreover, a continuum of elementary $(r,q)$-polycycles, which are proper
and convex on the hyperbolic plane.

The part (i) of Lemma above does not hold already for $(r,q)=(6,3)$: two elementary polycycles,
having each a single point as the kernel, can be glued into a $(6,3)$-polycycle
with $p_6=5$, having two isolated vertices as (disconnected) kernel, but
elementary polycycles, having those vertices as kernels, have a common 6-gon.

We believe that,
for each non-spheric $(r,q)$, amongst extremal $(r,q)$-polycycles 
there exist a proper one; in fact, each known extremal $(r,q)$-helicene 
$P$ has $p_r(P)>p_r(r^q)$.

The {\it extremal animals} of [\ref{HH}] are, in our terms, {\it proper} $(4,4$)-, 
$(6,3)$- and $(3,6)$-polycycles with minimal number of edges and so, maximal
number of interior vertices, for a given number of interior faces. It was proved
in [\ref{HH}] that such polycycles have $e=2p_4+ {\lceil 2 \sqrt{p_{4}} \rceil}$,
$e=3p_6+  {\lceil \sqrt{12p_{6}-3} \rceil}$, 
$e=p_3+ {\lceil \frac{p_3+\sqrt{6p_3}}{2} \rceil}$,
respectively, and that there are amongst them those, which grow like a spiral.
In fact, the proof of [\ref{HH}] implies that those extremal animals are, moreover, extremal in our sense.

Next two theorems give a full solution of the problem for $(r,q)=(5,3),(3,5)$.

\vspace{2mm}
\noindent
{\bf Theorem 7.}\ {\it Let $(r,q)=(5,3)$. Then

(i) with exception of 3 cases ($n(9)=10, n(10)=12, n(11)=15$),

$n(x)=x$, if $x \equiv 0,8,9 (mod \, 10)$, $n(x)=x-1$, if
$x \equiv 6,7 (mod \, 10)$, $n(x)=x-2$, otherwise,
and extremal polycycle is unique if $n(x)=x$,

(ii) all possible densities of $(5,3)$-polycycles, beside exceptions for
$p_5=9,10,11$, form the segment $ [0,1] $.}

The proof uses above Lemma and the list (found by exhaustive search) of all
connected components of kernels of $(5,3)$-polycycles; clearly, 
$E_i, 1 \le i \le 5,$ are first members of the family of
elementary $(5,3)$-polycycles $E_i, i \ge 1,$ with $i+2$ pentagons and $i$
interior vertices. 
The kernels are subgraphs induced by all interior vertices of elementary $(5,3)$- and $(3,5)$-polycycles on Figures 5 and 6.
Any $(5,3)$-polycycle can be obtained by gluing of elementary polycycles by
{\it open} (i.e. such that both vertices have degree 2) edges.
Theorem 7 is obtained by consideration of all such gluings.

\vspace{2mm}
\noindent
{\bf Theorem 8.}\ {\it Let $(r,q)=(3,5)$. Then

(i) $n(x)= \lfloor \frac{x}{3} \rfloor$, if $x \equiv 0,1 (mod \, 18)$,
with exception of $n(18)=8, n(19)=9$;

otherwise, $n(x)= \lfloor \frac{x-2}{3} \rfloor$,
with exception of $n(x)= \lfloor \frac{x+1}{3} \rfloor$ for

 $x=10,12,13,14,28,30,31,32,33,35$, and $n(x)= \lfloor \frac{x+4}{3} \rfloor$ for
$x=15,16,17,34$.

(ii) All possible densities of $(3,5)$-polycycles, besides of exceptions given
in (i) above, form the segment $ [ 0, \frac{1}{3}] $;
all rational densities can be realized by finite unproper polycycles.

(iii) For $p_3 \le 19$ the following holds: $n(x)=0$ if $0 \le x \le 4$,
 $n(x)=1$ if $5 \le x \le 7$, $n(x)=2$ if $8 \le x \le 9$, $n(x)=3$ if $10 \le x \le 11$,
$n(x)=4$ if $12 \le x \le 13$, $n(14)=5$, $n(x)=6$ if $15 \le x \le 16$,
$n(17)=7, n(18)=8, n(19)=9$.

Moreover, for $p_3 \le 19$, all extremal polycycles are proper and the 
reciprocal of an extremal one is also extremal; extremal ones are unique,
except of the cases $p_3=9, 11$ (two polycycles) and $p_3=4, 7, 13, 16$ (three polycycles).}

We used for the proof the same strategy as for the proof of Theorem 7,
but the number of cases to consider is much larger in this case. Other difference with $(5,3)$-polycycles is that we should use now (in order to glue some
elementary $(3,5)$-polycycles) the elementary polycycle $d$ with empty kernel.
So we give only Figure 6 of all minimal $(3,5)$-polycycles with connected 
kernels; clearly, $e_i, 1 \le i \le 6,$ are first members of the family of
elementary $(3,5)$-polycycles $e_i, i \ge 1,$ with $3i+2$ triangles and $i$
interior vertices.

\vspace{2mm}
Remark that the kernels of $(5,3)$-polycycles $A_i, 1 \le i \le 5,$ of 
Figure 5 are
all non-trivial isometric $(5,3)$-polycycles; they are also all 
{\it circumscribed} $(5,3)$-polycycles, i.e. $r$-gons can be
added around the perimeter, such that they will form a simple circuit.
(They were found in [{\ref{CC}]; such $(6,3)$-polycycles 
are useful in Organic Chemistry.) All circumscribed $(3,5)$-polycycles are the 
kernels of polycycles $a_i, 1 \le i \le 5,$ and $c_i, 1 \le i \le 4$ of 
Figure 6. It turns out that all polycycles $G$ in one of Figures 5 and 6, 
admitting the dual graph $G^*$ in another one, are as follows: $A^*_1=a_5$, $A^*_2=c_3$,
$A^*_3=c_4$, $A^*_4=e_2$, $A^*_5=e_1$; $E^*_1=d$ and $a^*_1=A_3$, $a^*_2=A_4$,
$a^*_4=C_3$, $a^*_5=A_5$, $c^*_1=E_4$, $c^*_2=E_3$, $c^*_3=E_2$, $c^*_4=E_1$;
$e^*_1=D$. $A_6$ occurs also in Figure 1; its dual
is infinite $(3,4)$-polycycle there. 

\newpage

\begin{center}
\input tab1a.pic
\input tab1b.pic
\end{center}
\begin{center}
{Figure 5: elementary (5,3)-polycycles and their kernels}
\end{center}

\newpage

\vspace{6mm}
\begin{center}
\input tab2a.pic
\input tab2b.pic
\end{center}
\begin{center}
{Figure 6: elementary (3,5)-polycycles and their kernels}
\end{center}

\newpage
For {\it any} $(r,q)$, we have only the bounds given in Theorem 9 below; 
they are linear with respect to $p_r$.
 The lower bound is attained, for example,  by the star of $q$ 
$r$-gons with one common (interior) vertex.
For Euclidean $(p,q)=(6,3),(4,4),(3,6)$, the upper bound is attained
in the limit.

\vspace{2mm}
\noindent
{\bf Theorem 9.}\ {\it Any $(r,q)$-polycycle, such that each its $r$-gon 
contains an interior vertex, satisfy to
$ \frac{p_r}{q} \le n_{int} < \frac{rp_r}{q}$.}

In order to prove those bounds, consider complex $X(G)$ 
from Section 3.2 above. Divide each its face (a regular $r$-gon) into $r$ $4$-gons by
the lines from the center to the mid-points of edges and, denoting by $\sigma$ 
the area of $4$-gon. Using that any $r$-gon contains at least one interior vertex, observe that

$p_r \sigma \le n_{int}q \sigma < p_rr \sigma$.

\subsection{Non-extendible polycycles}

Consider now another notion of maximality, appropriate to polycycles.
Call an $(r,q)$-polycycle {\it non-extendible} if it is not a proper subgraph
of another $(r,q)$-polycycle.
Four examples of non-extendible polycycles, depicted below, are: 
vertex-split $(3^4)$ and vertex-split $(3^5)$ (defined in Section 3) and two 
infinite polycycles ($Z$-paths of quadrangles and of triangles; both appear 
also in Figure 1).

\begin{center}
\input tab3.pic
\end{center}

\noindent
{\bf Theorem 10.}\ {\it 
All non-extendible $(r,q)$-polycycles are $(r^q)$, 4 above examples, 
possibly (but we conjecture their nonexistence) some other finite $(3,5)$-polycycles and, for any
$(r,q) \neq (3,3),(3,4),(4,3)$, a continuum of infinite ones.}

For example, a continuum of non-extendible $(5,3)$-polycycles comes as all 
infinite (in both directions) aperiodic sequences of glued
polycycles $C_2$ from Figure 5, and its upside down
version, say, $C_2^{'}$. The same procedure, used for polycycles $b_2, e_6$ from Figure 6,
gives a continuum of non-extendible $(3,5)$-polycycles. The $(4,4)$-cactus
from Figure 2 is also non-extendible since all its vertices have degree 4.

The hardest part of the proof of above Theorem was to show that any finite non-extendible polycycle $G$ is spheric.
Using a discrete analog of Gauss-Bonnet formula
(see Theorem 1.8.2 on page 76 of [\ref{Al}]), applied to $K(G)$), we get that 
the curvature of $K(G)$ is positive: the curvature is
equal to excess of the sum of angles of geodesic $n$-gon, bounding the disc
of perimeter $n$,
over $(n-2) \pi$, i.e. the sum of angles of Euclidean $n$-gon.
Namely, using non-extendibility of finite polycycle $G$, we get an estimation
of above sum of angles, which imply positivity of the curvature.

\vspace{2mm}
Finally, consider $(r,q)$-polycycles, such that any interior point of any
interior face has degree 1 of the cell-homomorphism onto $(r^q)$. The number of such polycycles, which are not extendible without loosing this property, is
equal to 0 for $(p,q)=(3,3),(3,4)$ and equal to 1 for $(p,q)=(4,3)$ (it is 
$P_2 \times P_5$). This number is finite for $(p,q)=(5,3),(3,5)$ and infinite,
otherwise.   

\section{Metric properties of polycycles: embedding into $Q_m$, $Z_m$}

Call a
polycycle {\it embeddable} if the metric space of its vertices (with usual
shortest-path metric) is embeddable isometrically, up to a scale  $ \lambda $, into
a hypercube  $Q_m$ , $m < \infty$, or (if the graph is infinite) into a
cubic lattice $Z_m$, $m \le \infty$; see [\ref{DL}],
[\ref{CDG}], [\ref{DS4}] for details. We have the following embeddings of
$(r^q)$: $(3^3)$, $(4^3)$ is embeddable into $Q_3$; $(3^4) into $ $Q_4$;
$(3^5)$ into $Q_{6}$; $(5^3)$ into $Q_{10}$; $(4^4)$ into $Z_2$; 
$(3^6)$, $(6^3)$ into $Z_3$ and all others into
$Z_ \infty$. So, any isometric polycycle is
embeddable.

Examples of non-embeddable polycycles are: $(4^3)-e$, $(3^4)-e$, $(5^3)-e$, 
$(3^5)-e$,
vertex-split $(3^4)$, vertex-split $(3^5)$ and four polycycles, given
on the Figures in Theorem below.
Amongst above 10 polycycles only vertex-split $(3^4)$ and vertex-split $(3^5)$ are
helicenes; amongst remaining 8 proper polycycles only ones on the right-hand side of Figures in 
Theorem below, are induced ones (they are $E_4$ and $c_3$ of Figures 5 and 6, respectively).

\vspace{2mm}
\noindent
{\bf Theorem 11.}{\it 

(i) for $(r,q) \neq (5,3),(3,5)$,
there are exactly 3 non-embeddable polycycles: $(4^3)-e$, $(3^4)-e$ and 
the vertex-split $(3^4)$,

(ii) except of $(5^3)$ itself, any 
 $(5,3)$-polycycle is embeddable if and only if it does not contain, as an
{\it induced} subgraph, neither of the two proper polycycles with $p_5=6$, given below

\medskip 
\begin{center}
\epsfig{file=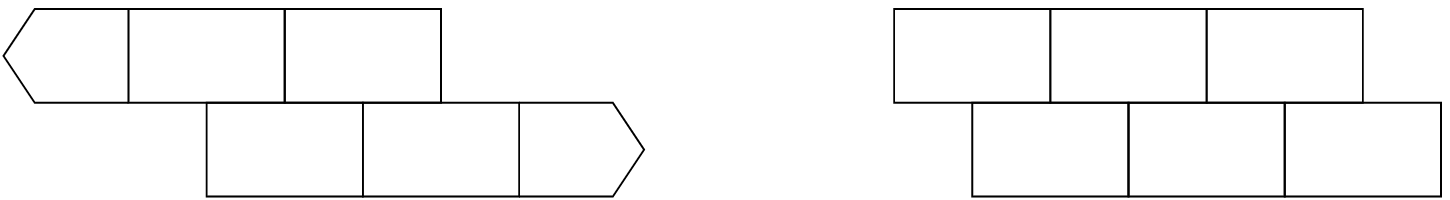,width=7.2cm}   
\end{center}
\medskip
  
(iii) we conjecture that, except of $(3^5)$ and $3^5-v$, any 
 $(3,5)$-polycycle is embeddable if and only if it does not contain, as an
{\it induced} subgraph, neither of the two proper 10-vertex polycycles, given 
below

\begin{center}
\input pl1.pic
\end{center}

There are no other non-embeddable $(3,5)$-polycycle with $p_3 \le 10$.}

\vspace{2mm}
Any outerplanar polycycle is embeddable, as well as any connected outerplanar graph.

Amongst 39 proper $(5,3)$-polycycles, the embeddable ones are: $(5^3)$, 3 with $p_5=7$, 9 
with $p_5=6$ (all but two given in (ii) above) and all 14 with $p_5 \le 5$.
Amongst $(5,3)$-polycycles of Figure 5, the embeddable ones are    
$A_1 \supset A_5$, $C_3 \supset E_3 \supset E_2 \supset E_1 \supset D$ 
(the sign ``$\supset$'' denotes here ``contains as a partial subgraph''). Amongst $(5,3)$-polycycles of 
Figure 1, only the infinite one with the symmetry $pma2$ and interior vertices, is not 
embeddable.

Amongst $(3,5)$-polycycles of Figure 6, the embeddable ones are $a_1 \supset a_5 \supset b_4 \supset c_4$, 
$e_3 \supset e_2 \supset e_1 \supset d$.
$(3,5)$-polycycles $a_5=(3^5)-v$, $a_2=(3^5)-e$ in Figure 6 contain 
the forbidden $(3,5)$-polycycles from (iii) above, as induced 
non-isometric subgraphs.

Embeddable $(r,q)$-polycycles have the scale $ \lambda = 1$ if
$r$ is even and $2$, otherwise. Consider any finite embeddable polycycle with perimeter $d$ and $k$
{\it closed zones}, i.e. cyclic sequences of opposite (alternating) non-boundary edges, 
see [\ref{CDG}], [\ref{DS4}]. Then it is embeddable into $Q_{ \frac{d}{2} + k}$, if $r$ is even (this implies that $d$ is
also even), and it is embeddable, with scale 2, into $Q_{d+k}$, if $r$ is odd. 
For example, amongst all embeddable polycycles in Figures 5 and 6, all with
$k>0$ are $A_1, a_1, a_5$; they have $k=5, 3, 1$, respectively.

There is a continuum of $(6,3)$-polycycles, which are embeddable only into
$Z_{ \infty}$: take all infinite (in both directions) paths of 
$(6,3)$-polycycles $P:=C_{0,1,...,9}+(3,8)$, glued each by its edges named 
$(4,5)$ and $(6,7)$; the choice how to glue - by edges of the same or
different name - is should be done aperiodically.

\section{Some relatives of $(r,q)$-polycycles}

It will be interesting to find some analog of above results 
for a generalization of polycycles on other surfaces. The following examples 
illustrate arising options:

(i) There is no straightforward analog of Theorem 1 for, say,
``polycycles on the torus $T^2$'': 2-connectedness of the graph does not imply
that its embedding on $T^2$ is a simply-connected union of faces (a handle can
broke it).  

(ii) The ring of three (or four) hexagons is example of planar, but not simply
connected ``polyhex'' (in a large sense used in Chemistry) without any (or locally-homeomorphic) homomorphism on 
$(6^3)$. Theorem 3 also does not admit a generalization already on {\it mono}-5-$(6^3)$, i.e.
on the partition of the plane by regular hexagons around one
central regular pentagon, having only vertices of degree 3 (in fact, a path of a
pentagon and six hexagons, having six
consecutive vertices of degree 3 on the boundary circuit, is not a
subgraph of  mono-5-$(6^3)$, but it is also not a ``helicene'', since, after 
mapping of its pentagon on the pentagon of  mono-5-$(6^3)$, the last hexagon
of the path also maps on this pentagon). On the other hand, the cell-map of Theorem 3 exists even
if we permit boundary vertices of degree greater than $q$, but this map
will be not locally-homeomorphic around those ``singularities''.

(iii) Consider {\it cross-ring}, depicted in Figure 8: first on $R^2$ as the
skeleton of dual disphenoid with 2 new vertices added on some 8 edges, and 
next in $R^3$ as
the skeleton of a polyhedral (2-connected) sphere with one hole and handle.
All 8 planar realizations of this graph are not polycyclic: they have both $8$-
and $9$-gonal faces. All conditions of the criterion in Theorem 1
are satisfied, except of (iv): $v-e+f=28-34+8 \neq 1$. But this polyhedral 
sphere permits to consider the cross-ring as ``not simply-connected 
$(8,3)$-polycycle'' with $p_8=5$; this polyhedron have five regular $8$-gonal 
faces: four indicated by the number ``8'' on the left-hand side of Figure 8 
and one bounded by $C_{11,12,4,5,25,26,18,19}$. The handle prevents it from 
embedding into $R^2$ and from cell-homomorphism into $(8^3)$.
All vertices $1,2,...,28$ are on the boundary of the disc.

\begin{center}
\epsfig{file=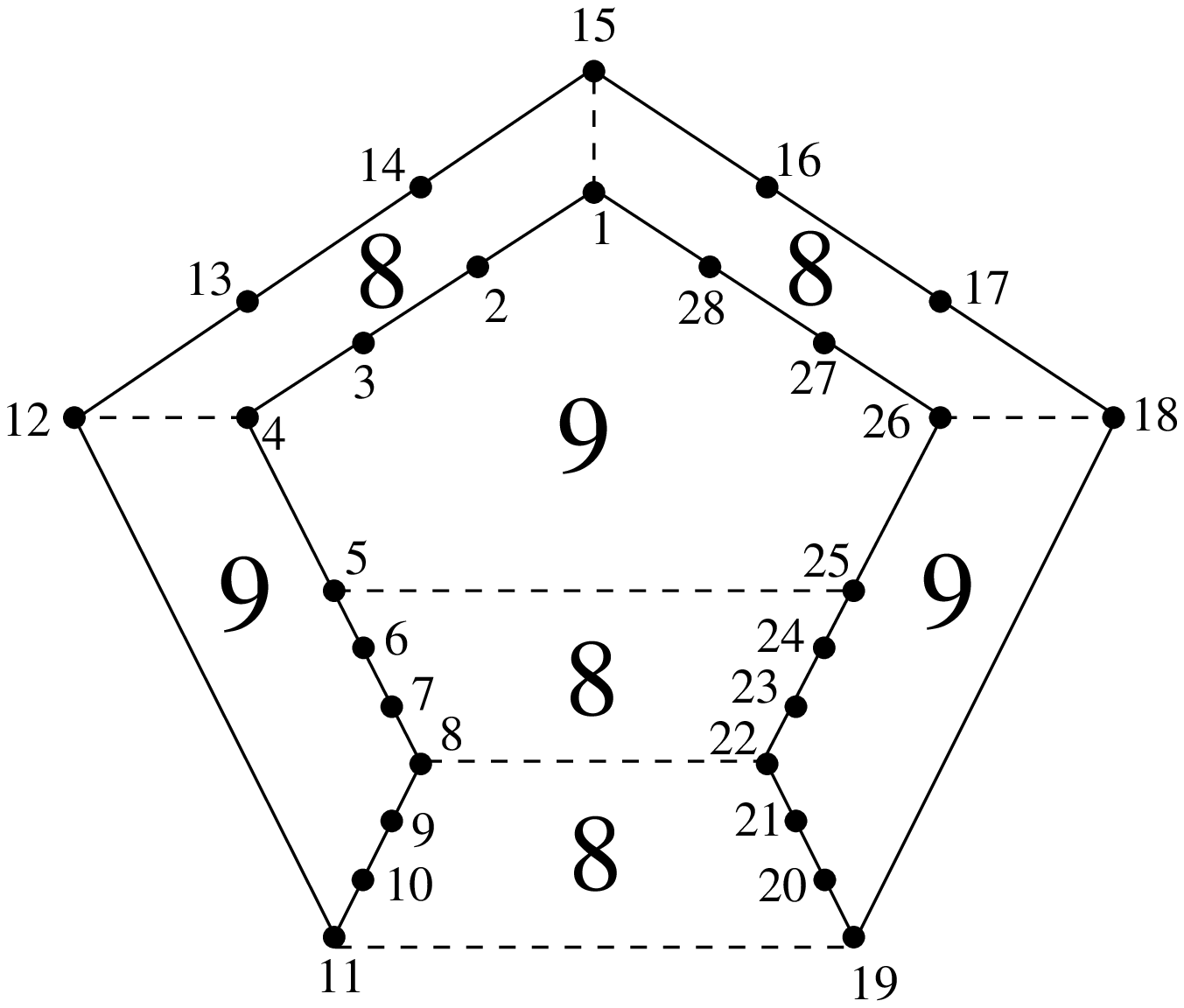,width=6.5cm}\hspace{2cm}\epsfig{file=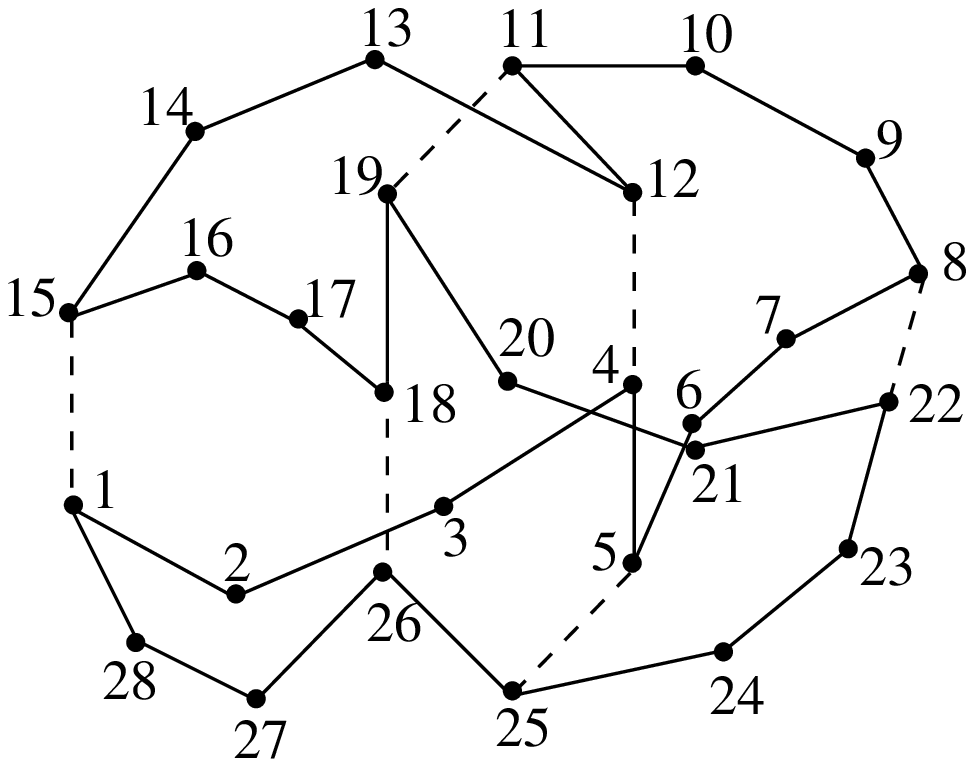,width=6.5cm}
\\[1cm]

\epsfig{file=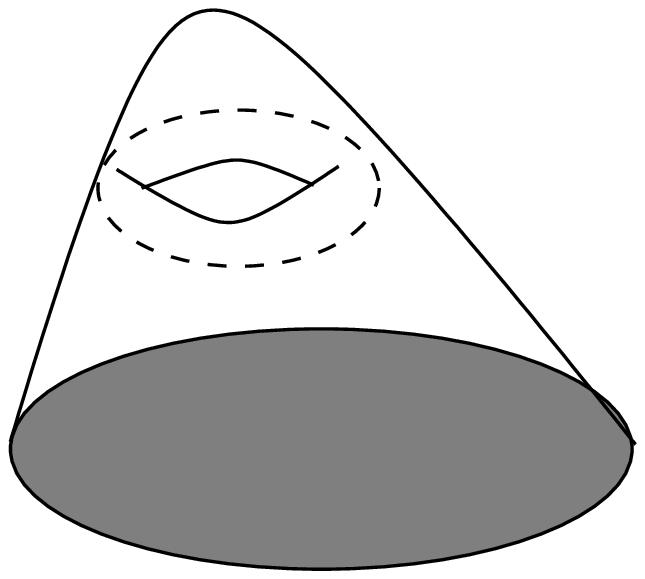,width=5cm} \\[1cm]

FIGURE 8: cross-ring as ``not simply-connected $(8,3)$-polycycle''
\end{center}

Finally, we mention two following relatives of finite $(r,q)$-polycycles.
See [\ref{AP}] and references there (mainly authored by M.Perkel) for the 
study of {\it strict polygonal} graphs, i.e. graphs of girth $r$ and
vertex-degree $q$, such that any path $P_2$ belongs to unique $r$-circuit 
of $G$. 
See pages 546-547 in [\ref{BW}] (and references there) for information on
{\it equivelar} polyhedra, i.e. polyhedral embedding, with convex faces,
of $(r,q)$-{\it map} (i.e. all faces are combinatorial $r$-gons and the 
vertex-degree of the skeleton is $q$) into $R^3$, such that all flags are equal.
So, in both these cases graphs are $q$-regular, have girth $r$,
Euler characteristic $v-e+f= \frac{v(6-r)}{2r}$ and coincide with Platonic polyhedra in the case of
genus 0. Recall that, for an $(r,q)$-polycycle, any non-boundary path 
$P_2$ belongs to unique $r$-circuit and there exist polyhedral realization
in $R^3$ with all interior faces being regular $r$-gons.

\newpage

\centerline{\Large References}
\vskip 10 pt
\begin{enumerate}
\small

\item
\label{Al}
A.D.Alexandrov (1950)
{\it Vypuklie mnogogranniki},
GITL, Moscow.

\item
\label{AP}
D.Archdeacon and M.Perkel (1990),
{\it Constructing Polygonal Graphs of Large Girth and Degree},
Congressus Numerantium, {\bf 70}, 81--85.

\item
\label{Ba}
 A.T.Balaban (1995),
 {\it Chemical Graphs: Looking back and Glimpsing ahead},
 J. Chem. Inf. Comput. Sci., {\bf 35}, 339--350.

\item
\label{BGOR}
 M.Bousquet-Mélou, A.J.Guttman, W.P.Orrick and A.Rechnitzer (1999),
{\it Inversion Relations, Reciprocity and Polyominoes}, Annals of Combinatorics, {\bf 3}, 223--249.

\item
\label{BW}
 U.Brehm and J.M.Wills (1993), {\it Polyhedral Manifolds}, Chapter 16 in
Handbook of Convex Geometry, edited by P.M.Gruber and J.M.Wills, Elsevier
Science Publishers.
 
\item
\label{CDG}
V.Chepoi, M.Deza and V.P.Grishukhin (1997)
{\it Clin d'oeil on $l_1$-embeddable
planar graphs}, Discrete Applied Math. {\bf 80}, 3--19.

\item
\label{CL}
 J.H.Conway and J.C.Lagarias (1990),
 {\it Tilings with Polyominoes and Combinatorial Group Theory},
 Journal of Combinatorial Theory, Series A  {\bf 53} (1990), 183--208.

\item
\label{CC}
 C.J.Cyvin, B.N.Cyvin, J.Brunvoll, E.Brendsdal, Zhang Fuji, Guo Xiofeng
 and R.Tosic (1993),
 {\it Theory of Polypentagons},
 J. Chem. Inf. Comput. Sci. {\bf 33}, 466--474.

\item
\label{DL}
 M.Deza and M.Laurent (1997),
 {\it Geometry of cuts and metrics},
 Springer-Verlag, Berlin.



\item
\label{DS3}
 M.Deza and M.I.Shtogrin (1998),
 {\it Polycycles},
 Voronoi Conference on Analytic Number Theory and Space Tilings (Kyiv,
September 7--14, 1998), Abstracts, Kyiv, 19--23.  

\item
\label{DS4}
 M.Deza and M.I.Shtogrin (2000),
 {\it Embedding of chemical graphs into hypercubes}, Math. Zametki {\bf 68-3}, 339--352.

\item
\label{DS5}
 M.Deza and M.I.Shtogrin (1999),
 {\it Primitive polycycles and helicenes},
 Russian Math. Surveys {\bf 54-6}, 159--160.

\item
\label{DS6}
 M.Deza and M.I.Shtogrin (2000),
 {\it Infinite primitive polycycles},
 Russian Math. Surveys {\bf 55-1}, 179--180.

\item
\label{DS7}
 M.Deza and M.I.Shtogrin (2000),
 {\it Polycycles: symmetry and embedding}, Russian Math. Surveys 
{\bf 56-6}, 159--160.

\item
\label{DS8}
 M.Deza and M.I.Shtogrin (2001),
 {\it Extremal and non-extendible polycycles}, Trudy Steklov Math. Institute, 
submitted.

\item
\label{Di}
 J.R.Dias (1988),
 {\it Handbook of polycyclic hydrocarbons. Part B: Polycyclic isomers
 and heteroatom analogs of benzenoid hydrocarbons}, Elsevier, Amsterdam.

\item
\label{GS}
B.Gr\"unbaum and G.C.Shephard (1986),
{\it Tilings and Patterns},
Freeman, New York.

\item
\label{Ha}
H.Harborth (1990),
{\it Some mosaic polyominoes}, ARS Combinatoria {\bf 29A}, 5--12.

\item
\label{HH}
F.Harary and H.Harborth (1976),
{\it Extremal animals}, Journal of Combinatorics, Information and System Sciences {\bf 1}, 1--8.

\item
\label{Ma}
W.Magnus (1974),
{\it Non-Euclidean Tessalations and their Groups},
Academic Press, New York and London.

\item
\label{S1}
M.I.Shtogrin (1999),
 {\it Primitive polycycles: criterion},
 Russian Math. Surveys {\bf 54-6}, 177--178.

\item
\label{S2}
M.I.Shtogrin (2000),
 {\it Non-primitive polycycles and helicenes},
 Russian Math. Surveys {\bf 55-2}, 155--156..

\end{enumerate}

{\it CNRS and Ecole Normale Sup\'erieure,
  45 rue d'Ulm, 75005 Paris, FRANCE. 

e-mail: Michel.Deza@ens.fr }

{\it Steklov Mathematical Institute,
8 Gubkin str., 117966 Moscow, RUSSIA.

e-mail: stogrin@mi.ras.ru }

\end{document}